\documentclass[12pt, a4paper]{article}
\usepackage{amsfonts,amsmath,amssymb,amsthm}
\usepackage{latexsym,amscd}
\usepackage{amsbsy}

\newcommand{\bcen}{\begin{center}}      \newcommand{\ecen}{\end{center}}
\def\B{{\cal B}}

\def\lz{\lambda}

\def\Hom{\mbox{Hom}}
\def\char{\mbox{char}}
\def\dim{\mbox{dim}}

\def\Ext{\mbox{Ext}}

\def\dim{\mbox{dim}}

\def\Im{\mbox{Im}}
\def\Ker{\mbox{Ker}}

\def\rank{\mbox{rank}}
\def\wt{\widetilde}

\def\wt{\widetilde}
\begin{document}

\begin{center}
{\Large {\bf Hochschild (co)homology of exterior
algebras\footnote{Project 10201004, 10426014 and 10501010
supported by NSFC}}}
\end{center}

\begin{center}
Yang Han $^{a}$ and Yunge Xu $^b$

\bigskip

{\footnotesize a. Academy of Mathematics and Systems science,
Chinese Academy of Sciences,\\ Beijing 100080, P.R. China
\hspace{5mm} E-mail: hany@iss.ac.cn

b. Faculty of Mathematics $\&$ Computer Science, Hubei University,\\
Wuhan 430062, P.R.China \hspace{5mm} E-mail: xuy@hubu.edu.cn}

\end{center}

\bigskip

\begin{center}

\begin{minipage}{12cm}
{\bf Abstract}: The minimal projective bimodule resolutions of the
exterior algebras are explicitly constructed. They are applied to
calculate the Hochschild (co)homology of the exterior algebras.
Thus the cyclic homology of the exterior algebras can be
calculated in case the underlying field is of characteristic zero.
Moreover, the Hochschild cohomology rings of the exterior algebras
are determined by generators and relations.

\medskip

{\bf Keyword:} Hochschild (co)homology, exterior algebra, minimal
projective resolution

\medskip

{\bf MSC(2000):} 16E40, 16G10

\end{minipage}

\end{center}

\section*{0. Introduction}

Fix a field $k$. Let $\Lambda$ be a finite-dimensional $k$-algebra
(associative with identity). Denote by $\Lambda^e$ the enveloping
algebra of $\Lambda$, i.e., the tensor product $\Lambda \otimes_k
\Lambda^{op}$ of the algebra $\Lambda$ and its opposite
$\Lambda^{op}$. The Hochschild homology and cohomology of
$\Lambda$ are defined by
$$HH_m(\Lambda) = \mbox{Tor}^{\Lambda^e}_{m}(\Lambda, \Lambda)\mbox{
and } HH^m(\Lambda)=\Ext^m_{\Lambda^e}(\Lambda, \Lambda)$$
respectively \cite{Mac}. The Hochschild (co)homology of an algebra
have played a fundamental role in representation theory of artin
algebras: Hochschild cohomology is closely related to simple
connectedness, separability and deformation theory
\cite{Skow,AP,Ger}; Hochschild homology is closely related to the
oriented cycle and the global dimension of algebras
\cite{I,LZ,AV,Han,Ke}.

Though Hochschild (co)homology is theoretically computable for a
concrete algebra via derived functors, actual calculation for a
class of algebras is still very convenient, important and
difficult. So far the Hochschild cohomology was calculated for
hereditary algebras \cite{C2,Hap}, incidence algebras
\cite{C1,GS}, algebras with narrow quivers \cite{Hap,C2}, radical
square zero algebras \cite{C3}, monomial algebras \cite{C4},
truncated quiver algebras \cite{C5,Zh1,Loc}, special biserial
algebras as well as their trivial extensions \cite{Xu,HX}, and so
on. The Hochschild homology was calculated for truncated algebras
\cite{LZ,Zh2,Skol}, quasi-hereditary algebras \cite{Za}, monomial
algebras \cite{Han}, and so on.

The exterior algebras play extremely important roles in many
mathematical branches such as algebraic geometry, commutative
algebra, differential geometry. It is well-known that the exterior
algebras can be viewed as $\mathbb{Z}/2$-graded algebras
\cite{Ka}. The graded exterior algebras (called {\it Grassmann
algebras} as well \cite{CJK}) have applications in physics and
their Hochschild (co)homology and cyclic (co)homology were known
(cf. \cite{Ka,CFRS,CJK,CR}). However, up to now the Hochschild
(co)homology of the ungraded exterior algebras (called {\it
parity-free Grassmann algebras} as well \cite{CJK}) are still
unknown. In this paper we shall deal with this problem. Our method
is purely algebraic and combinatorial. From now on all the
exterior algebras are ungraded. The content of this paper is
organized as follows: In section 1, we shall provide the minimal
projective bimodule resolutions of the exterior algebras. In
section 2, we shall apply these minimal projective bimodule
resolutions to calculate the Hochschild homology of the exterior
algebras and their cyclic homology in case the underlying field is
of characteristic zero. In section 3, we shall apply these minimal
projective bimodule resolutions to calculate the Hochschild
cohomology of the exterior algebras. In section 4, we shall
determine the Hochschild cohomology rings of the exterior algebras
by generators and relations.

\section{Minimal projective bimodule resolutions}

Let $Q$ be the quiver given by one point $1$ and $n$-loops
$x_1,x_2,...,x_n$ with $n \geq 2$. Denote by $I$ the ideal of the
path algebra $kQ$ generated by $R:=\{x_i^2|1 \leq i \leq n\} \cup
\{x_ix_j+x_jx_i| 1 \leq i<j \leq n\}$. For the knowledge on quiver
we refer to \cite{ARS}. Set $\Lambda=kQ/I$. Then $\Lambda$ is just
the exterior algebra over $k$ (cf. \cite{Mat}). Order the paths in
$Q$ by {\it left length lexicographic order} by choosing $1 < x_1
< x_2 < \cdots < x_n$, namely, $y_1 \cdots y_s < z_1 \cdots z_t$
with $y_i$ and $z_i$ being arrows if $s < t$ or if $s=t$, for some
$1 \leq r \leq s$, $y_i=z_i$ for $1 \leq i < r$ and $y_r<z_r$.
Then $\Lambda$ has a basis $\B=\cup_{i=0}^n\B_i$, where $\B_i=
\{x_{t_1}x_{t_2} \cdots x_{t_i}|1 \leq t_1 < t_2 < \cdots < t_i
\leq n \}$. So $\dim _k \Lambda = 2^n$. It is well-known that
$\Lambda$ is a Koszul algebra and its quadratic dual is just the
algebra of polynomials $k[x_1,...,x_n]$ (cf. \cite{BGS}).

Now we construct a minimal projective bimodule resolution
$(P_{\bullet},\delta_{\bullet})$ of $\Lambda$. Denote by $k
\langle x_1,...,x_n \rangle$ the noncommutative free associative
algebra over $k$ with free generators $x_1,...,x_n$. Denote by $k
\langle x_1,...,x_n \rangle _m$ the $k$-subspace of $k \langle
x_1,...,x_n \rangle$ generated by all monomials of degree $m$. For
each $m \geq 0$, we firstly construct elements $\{f^m_{1^{i_1}
2^{i_2}\cdots n^{i_n}} | i_1+i_2+\cdots +i_n=m, (i_1,...,i_n) \in
\mathbb{N}^n \} \subseteq k \langle x_1,...,x_n \rangle _m $: Let
$f_0^0=1, f^1_1=x_1, f^1_2=x_2, ..., f^1_n=x_n.$ Define
$f^m_{1^{i_1}2^{i_2}\cdots n^{i_n}}$ for all $m\geq 2$ inductively
by $f^m_{1^{i_1}2^{i_2}\cdots n^{i_{\scriptstyle
n}}}=\sum\limits_{h=1}^nf^{m-1}_{1^{i_1}\cdots h^{i_h-1}\cdots
n^{i_n}} x_h, $ where $i_1+i_2+\cdots +i_n=m, (i_1,...,i_n) \in
\mathbb{N}^n$ and $f^{m-1}_{1^{i_1}\cdots h^{-1}\cdots n^{i_n}}=0$
for all $1 \leq h \leq n$. It is well-known that the number of
non-negative integral solutions of the equation $i_1+i_2+\cdots
+i_n=m$ on $i_1,...,i_n$ is ${n+m-1\choose n-1}$ for any given
positive integers $n$ and $m$.

\medskip

Denote $\otimes := \otimes _k$. Let $P_m:=\coprod\limits_{i_1 +
i_2 + \cdots + i_n = m} \Lambda \otimes f^m_{1^{i_1}2^{i_2}\cdots
n^{i_n}} \otimes \Lambda \subseteq \Lambda \otimes k \langle
x_1,...,x_n \rangle _m \otimes \Lambda$ for $m \geq 0$, and let
$\wt{f}^m_{1^{i_1}2^{i_2}\cdots n^{i_{\scriptstyle n}}}:=1 \otimes
f^m_{1^{i_1}2^{i_2}\cdots n^{i_{\scriptstyle n}}} \otimes 1$ for
$m\geq 1$ and $\wt{f}^0_0=1\otimes 1$. Note that we identify $P_0$
with $\Lambda \otimes \Lambda$. Define $\delta_m: P_m\rightarrow
P_{m-1}$ by setting
$$
\delta_m(\wt{f}^m_{1^{i_1}2^{i_2}\cdots n^{i_{\scriptstyle n}}})=
\sum_{h=1}^n (x_h\wt{f}^{m-1}_{1^{i_1}\cdots h^{i_h-1}\cdots
n^{i_{\scriptstyle n}}} +(-1)^m\wt{f}^{m-1}_{1^{i_1}\cdots
h^{i_h-1}\cdots n^{i_{\scriptstyle n}}}x_h).
$$

\medskip

{\bf Theorem 1.} {\it The complex $(P_{\bullet},\delta_{\bullet})$
$$
\qquad \cdots\rightarrow P_{m+1}
\stackrel{\delta_{m+1}}{\longrightarrow} P_m
\stackrel{\delta_m}{\longrightarrow} \cdots
\stackrel{\delta_3}{\longrightarrow}
P_2\stackrel{\delta_2}{\longrightarrow} P_1
\stackrel{\delta_1}{\longrightarrow} P_0 \longrightarrow 0
$$
is a minimal projective bimodule resolution of the exterior
algebra $\Lambda=kQ/I$.}

\medskip

{\bf Proof.} For Koszul algebra $\Lambda$ we can construct its
Koszul complex which is a minimal projective resolution of the
trivial $\Lambda$-module $k$ (cf. \cite[Section 2.6]{BGS}).
Furthermore we can use the Koszul complex to construct the
bimodule Koszul complex which is a minimal projective bimodule
resolution of $\Lambda$ (cf. \cite[Section 9]{BK}). Let $X= k
\{x_1, x_2, \cdots, x_n\}$. We show that $\{f^m_{1^{i_1}
2^{i_2}\cdots n^{i_n}} | i_1+i_2+\cdots +i_n=m, (i_1,...,i_n) \in
\mathbb{N}^n \}$ is a $k$-basis of the $k$-vector space $K_m :=
\bigcap\limits_{p+q=m-2}X^pRX^q$ for all $m\geq 2$: Firstly, we
verify $f^m_{1^{i_1} 2^{i_2}\cdots n^{i_n}}\in K_m$: It is clear
that $f^m_{1^{i_1}2^{i_2}\cdots n^{i_n}}= \sum\limits_{h=1}^n
f^{m-1}_{1^{i_1} \cdots h^{i_h-1}\cdots n^{i_n}} x_h =
\sum\limits_{h=1}^n x_h f^{m-1}_{1^{i_1} \cdots h^{i_h-1} \cdots
n^{i_n}}. $ Thus the assertion follows by induction on $m$.
Secondly, $\{f^m_{1^{i_1} 2^{i_2}\cdots n^{i_n}} | i_1+i_2+\cdots
+i_n=m, (i_1,...,i_n) \in \mathbb{N}^n \}$ is $k$-linearly
independent: Indeed, by induction, one can show that each monomial
in $f^m_{1^{i_1} 2^{i_2}\cdots n^{i_n}}$  contains just $i_1$
$x_1$'s, $i_2$ $x_2$'s, ... , $i_n$ $x_n$'s. Thirdly, we have
$\dim_k K_m = {n+m-1 \choose n-1}$: The quadratic dual of the
Koszul algebra $\Lambda$ is just the algebra of polynomials
$k[x_1,...,x_n]$ which is isomorphic to the Yoneda algebra
$E(\Lambda)=\coprod\limits_{m \geq 0}\Ext^m_{\Lambda}(k,k)$ of
$\Lambda$ (cf. \cite[Theorem 2.10.1]{BGS}). Thus $\dim_k K_m =
\dim_k \Ext^m_{\Lambda}(k,k)= {n+m-1 \choose n-1}$. Hence
$\{f^m_{1^{i_1} 2^{i_2}\cdots n^{i_{\scriptstyle
n}}}|i_1+i_2+\cdots +i_n=m\}$ is a $k$-basis of $K_m$. Therefore
$P_{\bullet}$ are just those projective bimodules in the bimodule
Koszul complex of $\Lambda$ (cf. \cite[Section 9]{BK}).
Furthermore, $\delta_{\bullet}$ are just those differentials in
the bimodule Koszul complex of $\Lambda$ (cf. \cite[p. 354]{BK}).
\hfill$\square$

\section{Hochschild homology}

In this section we calculate the $k$-dimensions of Hochschild
homology groups and cyclic homology groups (in case char $k=0$) of
the exterior algebras.

\medskip

Applying the functor $\Lambda \otimes_{\Lambda^e}-$ to the minimal
projective bimodule resolution $(P_{\bullet},\delta_{\bullet})$,
we have $\Lambda
\otimes_{\Lambda^e}(P_{\bullet},\delta_{\bullet})=(M_{\bullet},
\tau_{\bullet})$ where $M_m=\coprod\limits_{i_1+i_2+\cdots+i_n=m}
\Lambda \otimes f^m_{1^{i_1} 2^{i_2}\cdots n^{i_{\scriptstyle
n}}}$, $\tau_m(\lambda \otimes f^m_{1^{i_1} 2^{i_2}\cdots
n^{i_n}})= \sum\limits_{h=1}^n (\lambda x_h \otimes
f^{m-1}_{1^{i_1} \cdots h^{i_h-1} \cdots n^{i_n}} +(-1)^m x_h
\lambda \otimes f^{m-1}_{1^{i_1} \cdots h^{i_h-1} \cdots
n^{i_n}})$ for $m \geq 1$, and $M_0=\Lambda \otimes k$. Throughout
we assume that the combinatorial number ${0 \choose 0}=1$ and ${i
\choose j}=0$ if $i < j$. Write $p(j)=p(m)$ if $j$ and $m$ are of
the same parity.

\medskip

{\bf Lemma 1.} {\it For} $m \geq 1$, {\it we have} $$\rank \,
\tau_m = \left\{\begin{array}{ll} \sum\limits^n_{i=1} {n \choose
i} \sum\limits^{i-1}_{j=0 \atop p(j)=p(m)} {j+m-1 \choose i-1}
{i-1 \choose j},&\mbox{if char $k \neq 2 $;}\\
0,&\mbox{otherwise.} \end{array}\right.$$

\medskip

{\bf Proof.} Denote by $k[x_1,...,x_n]_m$ the subspace of
$k[x_1,...,x_n]$ generated by all monomials of degree $m$.
Obviously the complex $(M_{\bullet}, \tau_{\bullet})$ is
isomorphic to the complex $(N_{\bullet},\sigma_{\bullet})$ which
is defined by $N_m:= \Lambda \otimes k[x_1,...,x_n]_m$ and
$\sigma_m : N_m \rightarrow N_{m-1}, \lambda \otimes x_1^{i_1}
\cdots x_n^{i_n} \mapsto \sum\limits_{h=1}^n (\lambda x_h \otimes
x_1^{i_1} \cdots x_h^{i_h-1} \cdots x_n^{i_n} +(-1)^m x_h \lambda
\otimes x_1^{i_1} \cdots x_h^{i_h-1} \cdots x_n^{i_n}).$

For any $\lambda=x_{t_1} \cdots x_{t_j} \in \B$, define
$\mu_{\lambda}(h):={|\{t_l | t_l < h, 1 \leq l \leq j \}|}$. Then
$$\sigma_m(\lambda \otimes x_1^{i_1} \cdots x_n^{i_n}) =
((-1)^j+(-1)^m) \sum\limits_{h=1}^n (-1)^{\mu_{\lambda}(h)}
N(\lambda x_h) \otimes x_1^{i_1} \cdots x_h^{i_h-1} \cdots
x_n^{i_n}$$ where $N(\lambda x_h):= x_{t_1} \cdots
x_{t_{\mu_{\lambda}(h)}}x_hx_{t_{\mu_{\lambda}(h)+1}} \cdots
x_{t_j}$, called the {\it normal form} of $\lambda x_h$ in
$\Lambda$. Note that the meaning of the normal form above is
completely different from that in the Gr\"{o}bner basis theory
(cf. \cite{G}). It follows from the formula above that
$\sigma_m(\lambda \otimes x_1^{i_1} \cdots x_n^{i_n}) =0$ if $p(j)
\neq p(m)$.

Clearly, $N_m$ has a basis ${\cal N}_m := \{ \lambda \otimes
x_1^{i_1} \cdots x_n^{i_n} \mid \lambda \in \B,\,
i_1+i_2+\cdots+i_n=m\}$. If $\lambda = x_{t_1} \cdots x_{t_j} \in
\B_j$ then $j$ is called the {\it degree} of $\lambda \otimes
x_1^{i_1} \cdots x_n^{i_n}$, and $j+ \sum\limits^n_{h=1}i_h$ is
called the {\it total degree} of $\lambda \otimes x_1^{i_1} \cdots
x_n^{i_n}$. If, viewed as a monomial in $k[x_1,...,x_n]$, $\lambda
x_1^{i_1} \cdots x_n^{i_n}$ can be written as $x_1^{j_1} \cdots
x_n^{j_n}$, then the set $\{ l | j_l \neq 0, 1 \leq l \leq n \} =
\{t_1, \cdots , t_j\} \cup \{ l | i_l \neq 0, 1 \leq l \leq n \}$
is called the {\it support} of $\lambda \otimes x_1^{i_1} \cdots
x_n^{i_n}$ and the cardinal number of the support $| \{ l | j_l
\neq 0, 1 \leq l \leq n \}|$ is called the {\it grade} of $\lambda
\otimes x_1^{i_1} \cdots x_n^{i_n}$.

Firstly, denote by $N_m(i)$ the subspace of $N_m$ generated by all
elements in ${\cal N}_m$ of grade $i$. Since $\sigma_m$ keeps the
grade of the elements, we have $\rank \, \sigma_m =
\sum\limits^n_{i=1} \rank (\sigma_m|N_m(i))$.

Secondly, denote by $N_m(x_{s_1} \cdots x_{s_i})$ the subspace of
$N_m(i)$ generated by all elements $\lambda \otimes x_1^{i_1}
\cdots x_n^{i_n}$ in ${\cal N}_m$ with support $\{s_1, \cdots ,
s_i \}$. Here $1 \leq s_1 < s_2 < \cdots < s_i \leq n$. The
restrictions of $\sigma _m$ to $N_m(x_{s_1} \cdots x_{s_i})$ and
$N_m(x_{s'_1} \cdots x_{s'_i})$ are given by the same matrices if
adequate bases are chosen, thus $\rank \, \sigma_m =
\sum\limits^n_{i=1} {n \choose i} \rank (\sigma_m|N_m(x_{1} \cdots
x_{i}))$.

Thirdly, denote by $N_m(x_1 \cdots x_i,j)$ the subspace of
$N_m(x_1 \cdots x_i)$ generated by all elements $\lambda \otimes
x_1^{l_1} \cdots x_i^{l_i}$ in ${\cal N}_m$ of degree $j$
(equivalently, of total degree $j+m$) with $0 \leq j \leq i$.
Moreover, denote by $N_m(x_1^{j_1} \cdots x_i^{j_i},j)$ the
subspace of $N_m(x_1 \cdots x_i,j)$ generated by all elements
$\lambda \otimes x_1^{l_1} \cdots x_i^{l_i}$ in ${\cal N}_m$ with
$\lambda x_1^{l_1} \cdots x_i^{l_i}= x_1^{j_1} \cdots x_i^{j_i}$,
$j_1,...,j_i \geq 1$ and $\sum\limits_{l=1}^ij_l=j+m$. Consider
the basis ${\cal N}_m(x_1^{j_1} \cdots x_i^{j_i},j) := {\cal N}_m
\cap N_m(x_1^{j_1} \cdots x_i^{j_i},j)$ of the vector space
$N_m(x_1^{j_1} \cdots x_i^{j_i},j)$. Order the elements $\lambda
\otimes x_1^{l_1} \cdots x_i^{l_i}$ in ${\cal N}_m(x_1^{j_1}
\cdots x_i^{j_i},j)$ by the left lexicographic order on $\lambda$.
Obviously $\sigma_m$ maps $N_m(x_1^{j_1} \cdots x_i^{j_i},j)$ into
$N_{m-1}(x_1^{j_1} \cdots x_i^{j_i},j+1)$ for $0 \leq j \leq i-1$
and into $0$ for $j=i$. Let $\alpha : N_m(x_1^{j_1} \cdots
x_i^{j_i},j) \rightarrow N_{m-1}(x_1^{j_1} \cdots x_i^{j_i},j+1)$
with $0 \leq j \leq i-1$ be the restriction of $\sigma_m$. Written
as a matrix under the basis ${\cal N}_m(x_1^{j_1} \cdots
x_i^{j_i},j)$ of $N_m(x_1^{j_1} \cdots x_i^{j_i},j)$ and ${\cal
N}_{m-1}(x_1^{j_1} \cdots x_i^{j_i},j+1)$ of $N_{m-1}(x_1^{j_1}
\cdots x_i^{j_i},j+1)$, $\alpha$ is a ${i \choose j+1} \times {i
\choose j}$ matrix. Partition the elements $\lambda \otimes
x_1^{l_1} \cdots x_i^{l_i}$ in ${\cal N}_m(x_1^{j_1} \cdots
x_i^{j_i},j)$ and ${\cal N}_{m-1}(x_1^{j_1} \cdots x_i^{j_i},j+1)$
according to whether $\lambda$ contains $x_1$ or not. In this way,
neglected the sign $(-1)^j+(-1)^m$, this matrix is partitioned
into a $2 \times 2$ partitioned matrix {\tiny $\left[\begin{array}{cc} A&I\\
0&B \end{array}\right]$} where $A$ is an ${i-1 \choose j} \times
{i-1 \choose j-1}$ matrix, $B$ is an ${i-1 \choose j+1} \times
{i-1 \choose j}$ matrix and $I$ is an ${i-1 \choose j} \times {i-1
\choose j}$ identity matrix.

\medskip

{\bf Claim.} $BA=0$.

\medskip

{\it Proof of the Claim.} Take any row of $B$. Assume that it
corresponds to the element $x_{q_1} \cdots x_{q_{j+1}} \otimes
x_1^{j_1} \cdots x_{q_1}^{j_{q_1}-1} \cdots
x_{q_{j+1}}^{j_{q_{j+1}}-1} \cdots x_i^{j_i} \in {\cal
N}_{m-1}(x_1^{j_1} \cdots x_i^{j_i},j+1)$. Then, according to the
definition of $\sigma_m$, this row is just the vector whose
components corresponding to $x_{q_1} \cdots x_{q_j} \otimes
x_1^{j_1} \cdots x_{q_1}^{j_{q_1}-1} \cdots x_{q_j}^{j_{q_j}-1}
\cdots x_i^{j_i}, ... ,$ \linebreak $ x_{q_2} \cdots x_{q_{j+1}}
\otimes x_1^{j_1} \cdots x_{q_2}^{j_{q_2}-1} \cdots
x_{q_{j+1}}^{j_{q_{j+1}}-1} \cdots x_i^{j_i}$ are
$(-1)^j,...,(-1)^0$ respectively, and other components are 0. Take
any column of $A$. Assume that it corresponds to the element
$x_1x_{q_1} \cdots \hat{x}_{q_a} \cdots \hat{x}_{q_b} \cdots
x_{q_{j+1}} \otimes x_1^{j_1-1} \cdots x_{q_1}^{j_{q_1}-1} \cdots
x_{q_a}^{j_{q_a}} \cdots $ \linebreak $x_{q_b}^{j_{q_b}} \cdots
x_{q_{j+1}}^{j_{q_{j+1}}-1} \cdots x_i^{j_i} \in {\cal
N}_m(x_1^{j_1} \cdots x_i^{j_i},j)$. Here $\hat{x}$ means $x$ is
deleted. This column is just the vector whose components
corresponding to $x_1x_{q_1} \cdots x_{q_a} \cdots $ \linebreak
$\hat{x}_{q_b} \cdots x_{q_{j+1}} \otimes x_1^{j_1-1} \cdots
x_{q_1}^{j_{q_1}-1} \cdots x_{q_a}^{j_{q_a}-1} \cdots
x_{q_b}^{j_{q_b}} \cdots x_{q_{j+1}}^{j_{q_{j+1}}-1} \cdots
x_i^{j_i}$ and $x_1x_{q_1} \cdots $ \linebreak $\hat{x}_{q_a}
\cdots x_{q_b} \cdots x_{q_{j+1}} \otimes x_1^{j_1-1} \cdots
x_{q_1}^{j_{q_1}-1} \cdots x_{q_a}^{j_{q_a}} \cdots
x_{q_b}^{j_{q_b}-1} \cdots x_{q_{j+1}}^{j_{q_{j+1}}-1} \cdots
x_i^{j_i}$ are $(-1)^a$ and $(-1)^{b-1}$ respectively, and other
components are 0. Thus the inner product of this row of $B$ and
this column of $A$ is just $(-1)^{b-1}(-1)^a +
(-1)^{a-1}(-1)^{b-1}=0$. Hence the claim holds.

\medskip

Note that under the chosen basis $\alpha$ is the matrix
$((-1)^j+(-1)^m)${\tiny $\left[\begin{array}{cc} A&I\\
0&B \end{array}\right]$}. Multiply {\tiny $\left[\begin{array}{cc} A&I\\
0&B \end{array}\right]$} on the left by the invertible matrix {\tiny $\left[\begin{array}{cc} I&0\\
-B&I \end{array}\right]$}, we obtain {\tiny $\left[\begin{array}{cc} I&0\\
-B&I \end{array}\right] \left[\begin{array}{cc} A&I\\
0&B \end{array}\right]=\left[\begin{array}{cc} A&I\\
-BA&0 \end{array}\right]=\left[\begin{array}{cc} A&I\\
0&0 \end{array}\right]$}. Hence rank{\tiny $\left[\begin{array}{cc} A&I\\
0&B \end{array}\right]$}$ = {i-1 \choose j}$. Thus
$$\begin{array}{lll} \rank (\sigma_m|N_m(x_1^{j_1} \cdots
x_i^{j_i},j))
& =& \rank \, \alpha \\
&=& \left\{\begin{array}{ll}
{i-1 \choose j},& \mbox{if $p(j)=p(m)$ and char $k \neq 2 $;}\\
0,&\mbox{otherwise.} \end{array}\right. \end{array}$$ Since
$$\begin{array}{lll} \rank (\sigma_m|N_m(x_1 \cdots x_i)) & = & \sum\limits^{i-1}_{j=0}
\sum\limits_{{j_1+ \cdots +j_i=j+m} \atop {j_1,...,j_i \geq 1}}
\rank (\sigma_m|N_m(x_1^{j_1} \cdots x_i^{j_i},j))\\
& = & \left\{\begin{array}{ll} \sum\limits^{i-1}_{j=0 \atop p(j)=p(m)}
{j+m-1 \choose i-1}{i-1 \choose j},&\mbox{if char $k \neq 2 $;}\\
0,&\mbox{otherwise.} \end{array}\right. \end{array}$$ we have
$$\begin{array}{lll} \rank \, \sigma_m & = & \sum\limits^n_{i=1} {n \choose i} \rank (\sigma_m|N_m(x_1 \cdots x_i))\\
& = & \left\{\begin{array}{ll} \sum\limits^n_{i=1} {n \choose i}
\sum\limits^{i-1}_{j=0 \atop p(j)=p(m)} {j+m-1 \choose i-1}{i-1 \choose j},&\mbox{if char $k \neq 2 $;}\\
0,&\mbox{otherwise.} \end{array}\right. \end{array}$$ The lemma
follows from the fact that $\rank \tau_{m} =\rank \sigma_{m}$.
\hfill{$\Box$}

\medskip

{\bf Lemma 2.} $\sum\limits_{i=j+1}^{n}{n \choose i}{j+m-1 \choose
i-1}{i-1 \choose j}=\sum\limits_{i=1}^{n-j}{n-i \choose j}{m+n-i-1
\choose n-i}$.

\medskip

{\bf Proof.} For $0 \leq j \leq n-1$, $m \geq 1-j$ and $0 \leq r
\leq i-1$, define $S_{m,\,r}=\sum\limits_{i=j+1}^{n}{n \choose
i}{j+m-1 \choose i-1-r}{i-1 \choose j},$ and
$T_{m,\,r}=\sum\limits_{i=1}^{n-j}{n-i \choose j}{m+n-i-1 \choose
n-i-r}.$ We shall show that $S_{m,\,r}=T_{m,\,r}$ for all $0 \leq
r \leq n-1$ and $m \geq 1-i$.

Firstly, we have $S_{m,\,n-1}=T_{m,\,n-1}$ for all $m \geq 1-j$:
This follows from $S_{m,\,n-1}= \sum\limits_{i=j+1}^{n}{n \choose
i} {j+m-1 \choose i-n}{i-1 \choose j} = {n \choose n}{j+m-1
\choose 0}{n-1\choose j} = {n-1 \choose j}$ and
$T_{m,\,n-1}=\sum\limits_{i=1}^{n-j}{n-i \choose j}{m+n-i-1
\choose 1-i}={n-1 \choose j}.$

Secondly, we have $S_{1-j,\,r}=T_{1-j,\,r}$ for all $0 \leq r \leq
n-1$: Note that $S_{1-j,\,r}=\sum\limits_{i=j+1}^{n}{n \choose
i}{0 \choose i-1-r}{i-1 \choose j}$ and $T_{1-j,\,r} =
\sum\limits_{i=1}^{n-j}{n-i \choose j}{n-i-j \choose n-i-r}.$ If
$r<i$ then $S_{1-j,\,r}=0=T_{1-j,\,r}.$ If $r \geq i$ then
$S_{1-j,\,r} = {n \choose r+1}{r\choose j}$ and $T_{1-j,\,r}=
\sum\limits_{i=1}^{n-j}{n-i \choose r}{r \choose i} ={n \choose
r+1}{r \choose j}$.

Thirdly, we have $S_{m,\,r}=T_{m,\,r}$ for all $0 \leq r \leq n-1$
and $m \geq 1-j$: Since $S_{m,\,r} = S_{m-1,\,r+1} + S_{m-1,\,r}$
and $T_{m,\,r} = T_{m-1,\,r+1} + T_{m-1,\,r}$, we have $$S_{m,\,r}
= S_{m-1,\,r+1}+S_{m-2,\,r+1}+\cdots+S_{1-j,\,r+1}+S_{1-j,\,r}$$
and $$T_{m,\,r} = T_{m-1,\,r+1} + T_{m-2,\,r+1} + \cdots +
T_{1-j,\,r+1} + T_{1-j,\,r}.$$ Thus
$$\begin{array}{ll}S_{m,\,n-2}
&= S_{m-1,\,n-1}+S_{m-2,\,n-1} +\cdots+S_{1-j,\,n-1}+S_{1-j,\,n-2}\\
&=T_{m-1,\,n-1}+T_{m-2,\,n-1}+\cdots+T_{1-j,\,n-1}+T_{1-j,\,n-2}\\
&=T_{m,\,n-2} \end{array}$$ for all $m \geq 1-j$. Similarly, by
induction, we have $S_{m,\,r}=T_{m,\,r}$ for all $0 \leq r \leq
n-1$ and $m \geq 1-j$.

Fourthly and finally, Lemma 2 follows from $S_{m,\,0}=T_{m,\,0}$.
\hfill$\square$

\medskip

{\bf Lemma 3.} {\it If} char$\,k \neq 2$ {\it and $m \geq 1$ then}
$$\rank \, \tau_m = \sum\limits_{i=1}^{n-1} 2^{i-1}{m+i-1\choose
i}+ \left\{\begin{array}{ll} 1,&\mbox{if $m$ is even;}\\
0,&\mbox{if $m$ is odd.}\end{array}\right.$$

\medskip

{\bf Proof.} By Lemma 2, we have
$$\begin{aligned} \rank \, \tau_m =&
\textstyle\sum\limits^{n-1}_{j=0 \atop
p(j)=p(m)}\sum\limits_{i=j+1}^{n}{n \choose i}{j+m-1 \choose
i-1}{i-1 \choose j}\\ =&\textstyle\sum\limits^{n-1}_{j=0 \atop
p(j)=p(m)}\sum\limits_{i=1}^{n-j}{n-i \choose j}{m+n-i-1\choose
n-i}\\ =&\textstyle\sum\limits_{i=1}^{n}\sum\limits^{n-i}_{j=0
\atop p(j)=p(m)}{n-i \choose j}{m+n-i-1 \choose n-i}\\ =&
\textstyle\sum\limits_{i=1}^{n}{m+n-i-1\choose n-i}
\sum\limits^{n-i}_{j=0 \atop p(j)=p(m)}{n-i \choose j}\\
=& \textstyle\sum\limits_{i=1}^{n-1}2^{n-i-1}{m+n-i-1\choose n-i}+
\left\{\begin{array}{ll} 1,&\mbox{if $m$ is even;}\\ 0,&\mbox{if $m$ is odd.}\end{array}\right.\\
=& \textstyle\sum\limits_{i=1}^{n-1}2^{i-1}{m+i-1\choose i}+
\left\{\begin{array}{ll} 1,&\mbox{if $m$ is even;}\\ 0,&\mbox{if
$m$ is odd.}\end{array}\right.
\end{aligned}$$
In the last step we replace $n-i$ with $i$. \hfill{$\square$}

\medskip

{\bf Lemma 4.} {\it If} char$\,k \neq 2$ {\it and $m \geq 1$ then}
$$\rank \, \tau_m + \rank \, \tau_{m+1} = 2^{n-1}{m+n-1\choose
n-1}.$$

\medskip

{\bf Proof.} Denote $U_n^m := \sum\limits_{i=1}^{n-1}2^{i-1}{m+i-1
\choose i}$. Then $U^{m+1}_n = \sum\limits_{i=1}^{n-1}2^{i-1} {m+i
\choose i}$ and $2U^{m+1}_n = \sum\limits_{i=1}^{n-1}2^i{m+i
\choose i}$. So $-U^{m+1}_n = {m+1 \choose 1} +
\sum\limits_{i=2}^{n-1}2^{i-1}{m+i-1 \choose i} - 2^{n-1}{m+n-1
\choose n-1} = U_n^m-2^{n-1}{m+n-1 \choose n-1}+1$. Thus
$U_n^m+U_n^{m+1}=2^{n-1}{m+n-1 \choose n-1}-1$. By Lemma 3, we
have $\rank \, \tau_m = U_n^m+ \left\{\begin{array}{ll}
1,&\mbox{if $m$ is even;}\\ 0,&\mbox{if $m$ is
odd.}\end{array}\right.$ Thus $\rank \, \tau_m + \rank \,
\tau_{m+1} = U_n^m + \left\{\begin{array}{ll} 1,&\mbox{if $m$ is
even;}\\ 0,&\mbox{if $m$ is odd.}\end{array}\right. +
U_n^{m+1} + \left\{\begin{array}{ll} 1,&\mbox{if $m+1$ is even;}\\
0,&\mbox{if $m+1$ is odd.}\end{array}\right. = 2^{n-1}{m+n-1
\choose n-1}.$

\hfill$\square$

\medskip

Now we can calculate the Hochschild homology of the exterior
algebra:

\medskip

{\bf Theorem 2.} {\it Let} $\Lambda=kQ/I$ {\it be the exterior
algebra. Then}
$$\dim_k HH_m(\Lambda) = \left\{\begin{array}{ll}
2^n {n+m-1 \choose n-1}, & \mbox{if char $k = 2$.}\\
2^{n-1}+1,& \mbox{if $m=0$ and char $k \neq 2$;}\\
2^{n-1}{n+m-1 \choose n-1}, & \mbox{if $m \geq 1$ and char $k \ne
2$;} \end{array}\right.$$

\medskip

{\bf Proof.} {\it Case} char $k\neq 2$: If $m \geq 1$ then, by
Lemma 4,
$$\begin{aligned}
\dim_k HH_m(\Lambda)=&\dim_k \, \Ker \, \tau_m-\dim_k \, \Im \, \tau_{m+1}\nonumber\\
=& \dim_k P_m-\dim_k \, \Im \, \tau_m-\dim_k \, \Im \, \tau_{m+1}\nonumber\\
=&\textstyle 2^n{n+m-1 \choose n-1}-\rank \, \tau_m-\rank \, \tau_{m+1}\nonumber\\
=&\textstyle 2^{n-1}{n+m-1 \choose n-1}
\end{aligned}$$

If $m=0$ then $HH_0(\Lambda) = \Lambda/[\Lambda,\Lambda]$, where
$[\Lambda,\Lambda] = \{\lz_1\lz_2-\lz_2\lz_1 | \lz_1, \lz_2 \in
\Lambda \}$. Note that as a vector space $[\Lambda,\Lambda]$ is
generated by all commutators $x_{i_1} \cdots x_{i_s}x_{j_1} \cdots
x_{j_t} - x_{j_1} \cdots x_{j_t}x_{i_1} \cdots x_{i_s}=
(1+(-1)^{st+1}) x_{i_1} \cdots x_{i_s}x_{j_1} \cdots x_{j_t}$ for
all $x_{i_1} \cdots x_{i_s}, x_{j_1} \cdots x_{j_t} \in {\cal B}$.
If both $s$ and $t$ are odd then $x_{i_1} \cdots x_{i_s}x_{j_1}
\cdots x_{j_t} \in [\Lambda,\Lambda]$. If $r \geq 2$ is even then
$r=(r-1)+1$. Thus $x_1 \cdots x_r \in [\Lambda,\Lambda]$. This
implies that $[\Lambda,\Lambda]$ has a basis consisting of the
images of all paths of even length ($\geq 2$). Therefore
$\Lambda/[\Lambda,\Lambda]$ has a basis consisting of the images
of identity and all paths of odd length. Hence
$\dim_kHH_0(\Lambda) = 1+ \sum\limits_{i=1 \atop {\scriptsize
\mbox{$i$ odd}}}^n {n \choose i} = 2^{n-1}+1$.

\medskip

{\it Case} char $k=2$: In this case all maps in the complex
$(M_{\bullet},\tau_{\bullet})$ are zero. Thus $HH_m(\Lambda) \cong
\Lambda^{{n+m-1 \choose n-1}}$ and $\dim_kHH_m(\Lambda) =
2^n{n+m-1 \choose n-1}$ for $m \geq 1$. Clearly
$\dim_kHH_0(\Lambda) = \dim_k\Lambda = 2^n$. \hfill $\square$

\medskip

Denote by $HC_m(\Lambda)$ the $m$-th cyclic homology group of $A$
(cf. \cite{Lod}). Let $hc_m(\Lambda):= \dim_k HC_m(\Lambda)$ and
$hh_m(\Lambda):= \dim_k HH_m(\Lambda)$

\medskip

{\bf Corollary 1.} {\it Let} $\Lambda=kQ/I$ {\it be the exterior
algebra and} $\char \, k =0$. {\it Then} $hc_m(\Lambda)=
\sum\limits^m_{i=0}
(-1)^{m-i}2^{n-1}{n+i-1 \choose n-1} + \left\{\begin{array}{ll} 1,&\mbox{if $m$ is even;}\\
0,&\mbox{if $m$ is odd.} \end{array}\right.$

\medskip

{\bf Proof.} By \cite[Theorem 4.1.13]{Lod}, we have
$$(hc_m(\Lambda)-hc_m(k)) = -(hc_{m-1}(\Lambda)-hc_{m-1}(k))+
(hh_m(\Lambda)-hh_m(k)).$$ Thus $(hc_m(\Lambda)-hc_m(k)) =
\sum\limits^m_{i=0}(-1)^{m-i}(hh_i(\Lambda)-hh_i(k)).$ It is
well-known that $hh_i(k)=\left\{\begin{array}{ll} 1,&\mbox{if
$i=0$;}\\ 0,&\mbox{if $i \geq 1$.} \end{array}\right.$ and
$hc_m(k)=\left\{\begin{array}{ll} 1,&\mbox{if $m$ is even;}\\
0,&\mbox{if $m$ is odd.} \end{array}\right.$ By Theorem 2, we have
$hc_m(\Lambda)= \sum\limits^m_{i=0}
(-1)^{m-i}2^{n-1}{n+i-1 \choose n-1} + \left\{\begin{array}{ll} 1,&\mbox{if $m$ is even;}\\
0,&\mbox{if $m$ is odd.} \end{array}\right.$ \hfill{$\Box$}

\section{Hochschild cohomology}

In this section we calculate the $k$-dimensions of the Hochschild
cohomological groups of the exterior algebras.

Applying the functor $\Hom_{\Lambda^e}(-,\Lambda)$ to the minimal
projective bimodule resolution $(P_{\bullet},\delta_{\bullet})$,
we have $\Hom_{\Lambda^e}((P_{\bullet},\delta_{\bullet}), \Lambda)
=(P_{\bullet}^*, \delta_{\bullet}^*)$ where $P^*_m =
\Hom_{\Lambda^e}(P_m,\Lambda)$ and $\delta^*_m(\phi) = \phi
\delta_m$ for any $\phi \in P^*_m$. As $k$-vector spaces, $P^*_m$
is canonically isomorphic to $M^m := \Lambda^{{n+m-1 \choose
n-1}}$. Let $\{\phi ^m_{1^{i_1} \cdots n^{i_n}} | i_1+ \cdots +i_n
= m\}$ be a basis of the $k$-vector space $P^*_m$ defined by $\phi
^m_{1^{i_1} \cdots n^{i_n}}(1 \otimes f^m_{1^{i_1} \cdots n^{i_n}}
\otimes 1)=1$ and $\phi ^m_{1^{i_1} \cdots n^{i_n}}(1 \otimes
f^m_{1^{j_1} \cdots n^{j_n}} \otimes 1)=0$ for $(j_1, \cdots ,
j_n) \neq (i_1, \cdots , i_n)$. Let $e^m_{1^{i_1} \cdots n^{i_n}}$
be the image of $\phi ^m_{1^{i_1} \cdots n^{i_n}}$ under the
canonical isomorphism $P^*_m \cong M^m$. Then the complex
$(P_{\bullet}^*, \delta_{\bullet}^*)$ is isomorphic to the complex
$(M^{\bullet},\tau^{\bullet})$ where $\tau^{m+1}: M^m \rightarrow
M^{m+1}, \lambda e^m_{1^{i_1} \cdots n^{i_n}} \mapsto
\sum\limits^n_{h=1}(x_h\lambda + (-1)^{m+1} \lambda x_h)
e^m_{1^{i_1} \cdots h^{i_h+1} \cdots n^{i_n}} =
(1+(-1)^{m+j+1})\sum\limits^n_{h=1}(-1)^{\mu_{\lambda}(h)}N(x_h\lambda)
e^{m+1}_{1^{i_1} \cdots h^{i_h+1} \cdots n^{i_n}}$ for $\lambda
\in {\cal B}_j$.

\medskip

{\bf Lemma 5.} {\it For} $m \geq 0$ {\it we have}
$$\begin{array}{lll} \rank \, \tau^{m+1} \!\! & \!\! = \!\! & \!\!
\left\{\begin{array}{ll} \sum\limits^n_{i=1} {n \choose i}
\sum\limits^{i-1}_{j=0 \atop
p(j)=p(n+m)} {j+m \choose i-1} {i-1 \choose j},&\mbox{if char $k \neq 2 $;}\\
0,&\mbox{otherwise.} \end{array}\right. \end{array}$$

\medskip

{\bf Proof.} Obviously the complex $(M^{\bullet}, \tau^{\bullet})$
is isomorphic to the complex $(N^{\bullet},\sigma^{\bullet})$
which is defined by $N^m:= \Lambda \otimes k[x_1,...,x_n]_m$ and
$\sigma^{m+1} : N^m \rightarrow N^{m+1}, \lambda \otimes x_1^{i_1}
\cdots x_n^{i_n} \mapsto \sum\limits_{h=1}^n (x_h \lambda \otimes
x_1^{i_1} \cdots x_h^{i_h+1} \cdots x_n^{i_n} +(-1)^{m+1} \lambda
x_h \otimes$ \linebreak $x_1^{i_1} \cdots x_h^{i_h+1} \cdots
x_n^{i_n}).$ For any $\lambda=x_{t_1} \cdots x_{t_j} \in \B$, we
have $\sigma^{m+1}(\lambda \otimes x_1^{i_1} \cdots x_n^{i_n})$ $
= (1+(-1)^{m+j+1}) \sum\limits_{h=1}^n (-1)^{\mu_{\lambda}(h)}
N(\lambda x_h) \otimes x_1^{i_1} \cdots x_h^{i_h+1} \cdots
x_n^{i_n}$.

Clearly, $N^m$ has a basis ${\cal N}^m := \{ \lambda \otimes
x_1^{i_1} \cdots x_n^{i_n} \mid \lambda \in \B,\,
i_1+i_2+\cdots+i_n=m\}$. If, viewed as a monomial in
$k[x_1,...,x_n,x_1^{-1},...,x_n^{-1}]$, $\lambda x_1^{-i_1} \cdots
x_n^{-i_n}$ can be written as $x_1^{j_1} \cdots x_n^{j_n}$, then
the sum $i$ (resp. $-j$) of all the positive (resp. negative)
$j_l$ is called the {\it positive} (resp. {\it negative}) {\it
degree} of $\lambda \otimes x_1^{i_1} \cdots x_n^{i_n}$. Note that
the positive $j_l$ must be 1. Moreover, $\{l | j_l > 0, 1 \leq l
\leq n \}$ (resp. $\{l | j_l < 0, 1 \leq l \leq n \}$) is called
the {\it positive} (resp. {\it negative}) {\it support} of
$\lambda \otimes x_1^{i_1} \cdots x_n^{i_n}$.

Firstly, denote by $N^m(i)$ the subspace of $N^m$ generated by all
elements in ${\cal N}^m$ of positive degree $i$. Since
$\sigma^{m+1}$ keeps the positive degree of the elements, we have
$\rank \, \sigma^{m+1} = \sum\limits^n_{i=1} \rank
(\sigma^{m+1}|N^m(i))$.

Secondly, denote by $N^m(x_{s_1} \cdots x_{s_i})$ the subspace of
$N^m(i)$ generated by all elements $\lambda \otimes x_1^{i_1}
\cdots x_n^{i_n}$ in ${\cal N}^m$ with positive support $\{s_1,
\cdots , s_i \}$. Here $1 \leq s_1 < s_2 < \cdots < s_i \leq n$.

Thirdly, denote by $N^m(x_{s_1} \cdots
x_{s_i},x_{s_{i+1}}^{j_{i+1}} \cdots x_{s_n}^{j_n})$, where $0
\leq i \leq n, 0 \leq j \leq m, j_{i+1}+ \cdots +j_n=j,
\{s_{i+1},...,s_n\}=\{1,...,n\} \backslash \{s_1,...,s_i\}$ and
$s_{i+1}< \cdots <s_n$, the subspace of $N^m$ generated by all
elements $N(x_{s_1} \cdots x_{s_i}x_{s_{t_1}} \cdots
x_{s_{t_{m-j}}}) \otimes x_{s_{t_1}} \cdots
x_{s_{t_{m-j}}}x_{s_{i+1}}^{j_{i+1}} \cdots x_{s_n}^{j_n}$ of
positive degree $i$ and negative degree $j$ in ${\cal N}^m$ with
$\{s_{t_1}, ... ,s_{t_{m-j}}\} \subseteq \{s_{i+1},...,s_{n}\}$
and $s_{t_1} < \cdots <s_{t_{m-j}}$. Clearly this subspace is 0
unless $0 \leq m-j \leq n-i$. Consider the basis ${\cal
N}^m(x_{s_1} \cdots x_{s_i},$ \linebreak $x_{s_{i+1}}^{j_{i+1}}
\cdots x_{s_n}^{j_n}) := {\cal N}^m \cap N^m(x_{s_1} \cdots
x_{s_i},x_{s_{i+1}}^{j_{i+1}} \cdots x_{s_n}^{j_n})$ of the vector
space \linebreak $N^m(x_{s_1} \cdots x_{s_i},x_{s_{i+1}}^{j_{i+1}}
\cdots x_{s_n}^{j_n})$. Order the elements $\lambda \otimes
x_{s_{i+1}}^{l_{i+1}} \cdots x_{s_n}^{l_n}$ in \linebreak ${\cal
N}^m(x_{s_1} \cdots x_{s_i},x_{s_{i+1}}^{j_{i+1}} \cdots
x_{s_n}^{j_n})$ by the left lexicographic order on $\lambda$.
Obviously $\sigma^{m+1}$ maps $N^m(x_{s_1} \cdots
x_{s_i},x_{s_{i+1}}^{j_{i+1}} \cdots x_{s_n}^{j_n})$ into
$N^{m+1}(x_{s_1} \cdots x_{s_i},x_{s_{i+1}}^{j_{i+1}} \cdots
x_{s_n}^{j_n})$ if $0 \leq i \leq n-1$, and into 0 if $i=n$. Let
the map $\alpha : N^m(x_{s_1} \cdots x_{s_i},x_{s_{i+1}}^{j_{i+1}}
\cdots x_{s_n}^{j_n})$ \linebreak $ \rightarrow N^{m+1}(x_{s_1}
\cdots x_{s_i},x_{s_{i+1}}^{j_{i+1}} \cdots x_{s_n}^{j_n})$ be the
restriction of $\sigma^{m+1}$. Written as a matrix under the basis
${\cal N}^m(x_{s_1} \cdots x_{s_i},x_{s_{i+1}}^{j_{i+1}} \cdots
x_{s_n}^{j_n})$ of $N^m(x_{s_1} \cdots
x_{s_i},x_{s_{i+1}}^{j_{i+1}} \cdots x_{s_n}^{j_n})$ and ${\cal
N}^{m+1}(x_{s_1} \cdots x_{s_i},x_{s_{i+1}}^{j_{i+1}} \cdots
x_{s_n}^{j_n})$ of $N^{m+1}(x_{s_1} \cdots
x_{s_i},x_{s_{i+1}}^{j_{i+1}} \cdots x_{s_n}^{j_n})$, $\alpha$ is
an ${n-i \choose m-j+1} \times {n-i \choose m-j}$ matrix.
Partition the elements $\lambda \otimes x_{s_{i+1}}^{l_{i+1}}
\cdots x_{s_n}^{l_n}$ in \linebreak ${\cal N}^m(x_{s_1} \cdots
x_{s_i},x_{s_{i+1}}^{j_{i+1}} \cdots x_{s_n}^{j_n})$ and ${\cal
N}^{m+1}(x_{s_1} \cdots x_{s_i},x_{s_{i+1}}^{j_{i+1}} \cdots
x_{s_n}^{j_n})$ according to whether $\lambda$ contains
$x_{s_{i+1}}$ or not. In this way, neglected the sign \linebreak
$1+(-1)^{m+(i+m-j)+1}=1+(-1)^{i-j+1}$, this matrix is partitioned
into a $2 \times 2$ partitioned matrix
{\tiny $\left[\begin{array}{cr} A&(-1)^{\mu_{x_{s_1} \cdots x_{s_i}}(s_{i+1})}I\\
0&B \end{array}\right]$} where $A$ is an ${n-i-1 \choose m-j}
\times {n-i-1 \choose m-j-1}$ matrix, $B$ is an ${n-i-1 \choose
m-j+1} \times {n-i-1 \choose m-j}$ matrix and $I$ is an ${n-i-1
\choose m-j} \times {n-i-1 \choose m-j}$ identity matrix.

\medskip

{\bf Claim.} $BA=0$.

\medskip

{\it Proof of the Claim.} Take any row of $B$. Assume that it
corresponds to the element $N(x_{s_1} \cdots x_{s_i}x_{s_{q_1}}
\cdots x_{s_{q_{j+1}}}) \otimes x_{s_{i+1}}^{j_{i+1}} \cdots
x_{s_{q_1}}^{j_{q_1}+1} \cdots x_{s_{q_{j+1}}}^{j_{q_{j+1}}+1}
\cdots x_{s_n}^{j_n} \in {\cal N}^{m+1}(x_{s_1} \cdots
x_{s_i},x_{s_{i+1}}^{j_{i+1}} \cdots x_{s_n}^{j_n})$. Then,
according to the definition of $\sigma^{m+1}$, this row is just
the vector whose components corresponding to $N(x_{s_1} \cdots
x_{s_i}x_{s_{q_1}} \cdots $ \linebreak $x_{s_{q_j}}) \otimes
x_{s_{i+1}}^{j_{i+1}} \cdots x_{s_{q_1}}^{j_{q_1}+1} \cdots
x_{s_{q_j}}^{j_{q_j}+1} \cdots x_{s_n}^{j_n}$, ... , $N(x_{s_1}
\cdots x_{s_i}x_{s_{q_2}} \cdots x_{s_{q_{j+1}}}) \otimes$
\linebreak $ x_{s_{i+1}}^{j_{i+1}} \cdots x_{s_{q_2}}^{j_{q_2}+1}
\cdots x_{s_{q_{j+1}}}^{j_{q_{j+1}}+1} \cdots x_{s_n}^{j_n}$ are
$(-1)^{\mu_{x_{s_1} \cdots
x_{s_i}}(s_{q_{j+1}})+j},...,(-1)^{\mu_{x_{s_1} \cdots
x_{s_i}}(s_{q_1})}$ \linebreak respectively, and other components
are 0. Take any column of $A$. Assume that it corresponds to the
element $N(x_{s_1} \cdots x_{s_{i+1}}x_{s_{q_1}} \cdots
\hat{x}_{s_{q_a}} \cdots \hat{x}_{s_{q_b}} \cdots
x_{s_{q_{j+1}}})$ \linebreak $ \otimes x_{s_{i+1}}^{j_{i+1}+1}
\cdots x_{s_{q_1}}^{j_{q_1}+1} \cdots x_{s_{q_a}}^{j_{q_a}} \cdots
x_{s_{q_b}}^{j_{q_b}} \cdots x_{s_{q_{j+1}}}^{j_{q_{j+1}}+1}
\cdots x_{s_n}^{j_n} \in {\cal N}^m(x_{s_1} \cdots x_{s_i},
x_{s_{i+1}}^{j_{i+1}} \cdots $ \linebreak $x_{s_n}^{j_n})$. This
column is just the vector whose components corresponding to
$N(x_{s_1} \cdots x_{s_{i+1}}x_{s_{q_1}} \cdots x_{s_{q_a}} \cdots
\hat{x}_{s_{q_b}} \cdots x_{s_{q_{j+1}}}) \otimes
x_{s_{i+1}}^{j_{i+1}+1} \cdots x_{s_{q_1}}^{j_{q_1}+1} \cdots
x_{s_{q_a}}^{j_{q_a}+1} \cdots $ \linebreak $x_{s_{q_b}}^{j_{q_b}}
\cdots x_{s_{q_{j+1}}}^{j_{q_{j+1}}+1} \cdots x_{s_n}^{j_n}$ and
$N(x_{s_1} \cdots x_{s_{i+1}}x_{s_{q_1}} \cdots \hat{x}_{s_{q_a}}
\cdots x_{s_{q_b}} \cdots x_{s_{q_{j+1}}}) \otimes$ \linebreak $
x_{s_{i+1}}^{j_{i+1}+1} \cdots x_{s_{q_1}}^{j_{q_1}+1} \cdots
x_{s_{q_a}}^{j_{q_a}} \cdots x_{s_{q_b}}^{j_{q_b}+1} \cdots
x_{s_{q_{j+1}}}^{j_{q_{j+1}}+1} \cdots x_{s_n}^{j_n}$ are
$(-1)^{\mu_{x_{s_1} \cdots x_{s_i}}(s_{q_a})+a}$ and
$(-1)^{\mu_{x_{s_1} \cdots x_{s_i}}(s_{q_b})+b-1}$ respectively,
and other components are 0. Thus the inner product of this row of
$B$ and this column of $A$ is just \linebreak $(-1)^{\mu_{x_{s_1}
\cdots x_{s_i}}(s_{q_b})+b-1}(-1)^{\mu_{x_{s_1} \cdots
x_{s_i}}(s_{q_a})+a} + (-1)^{\mu_{x_{s_1} \cdots
x_{s_i}}(s_{q_a})+a-1}(-1)^{\mu_{x_{s_1} \cdots
x_{s_i}}(s_{q_b})+b-1}$ \linebreak $=0$. Hence the claim holds.

\medskip

Under the chosen basis, $\alpha$ is the matrix
$(1+(-1)^{i-j+1})${\tiny $\left[\begin{array}{cr}
A&(-1)^{\mu_{x_{s_1} \cdots x_{s_i}}(s_{i+1})}I\\
0&B \end{array}\right]$}. Thus we have $\rank \, \alpha =
\left\{\begin{array}{ll}
{n-i-1 \choose m-j},& \mbox{if $p(i) \ne p(j)$ and char $k \neq 2 $;}\\
0,&\mbox{otherwise.} \end{array}\right.$ Since $\sigma^{m+1}$
preserves both positive degree and negative degree, we have

$$\begin{array}{ll} &\rank \, \sigma^{m+1}\\ =& \sum\limits^n_{i=0} \rank
(\sigma^{m+1}|N^m(i))\\ =&  \sum\limits^n_{i=0}
\sum\limits_{\{s_1,...,s_i\} \subseteq \{1,...,n\}} \rank
(\sigma^{m+1}|N^m(x_{s_1} \cdots x_{s_i}))\\  =&
\sum\limits^n_{i=0} \sum\limits_{\{s_1,...,s_i\} \subseteq
\{1,...,n\}} \sum\limits^m_{j=m-n+i}  \sum\limits_{j_{i+1}+ \cdots
+j_n=j} \rank (\sigma^{m+1}|N^m(x_{s_1} \cdots
x_{s_i},x_{s_{i+1}}^{j_{i+1}} \cdots x_{s_n}^{j_n}))\\
 =& \left\{\begin{array}{ll} \sum\limits^{n-1}_{i=0}
{n \choose i} \sum\limits^m_{j=m-n+i \atop p(j) \ne p(i)} {n-i+j-1
\choose
n-i-1}{n-i-1 \choose m-j},&\mbox{if char $k \neq 2 $;}\\
0,&\mbox{otherwise.} \end{array}\right. \end{array}$$ If char $k
\neq 2$ then we have $$\begin{array}{lll} \rank \sigma^{m+1}& =&
\sum \limits^{n-1}_{i=0} {n \choose i} \sum\limits^m_{j=m-n+i
\atop p(j) \ne p(i)} {n-i+j-1 \choose n-i-1} {n-i-1 \choose m-j}\\
& = & \sum \limits^n_{i=1} {n \choose i} \sum\limits^i_{j=0 \atop
p(m-j) \ne p(n-i)} {m-j+i-1 \choose i-1} {i-1 \choose j}\\ & = &
\sum \limits^n_{i=1} {n \choose i} \sum\limits^i_{j=0 \atop p(j)
\ne p(n+m)} {m+j-1 \choose i-1} {i-1 \choose j-1}
\\ & = & \sum \limits^n_{i=1} {n \choose i} \sum\limits^{i-1}_{j=0
\atop p(j) = p(n+m)} {m+j \choose i-1} {i-1 \choose j}.
\end{array}$$ The lemma follows from the fact that $\rank \, \tau^{m+1}
=\rank \, \sigma^{m+1}$. \hfill{$\Box$}

\medskip

{\bf Lemma 6.}  {\it If} char$\,k \neq 2$ {\it and $m \geq 1$
then}
$$\rank \, \tau^{m}+\rank \, \tau^{m+1}=2^{n-1}{n+m-1 \choose
n-1}.$$

\medskip

{\bf Proof.} Note that
$$\begin{aligned} \rank \, \tau^{m+1}=& \textstyle\sum\limits^n_{i=1} {n \choose i}
\sum\limits^{i-1}_{j=0 \atop p(j) = p(n+m)} {j+m \choose i-1} {i-1
\choose j}\\ = & \textstyle\sum\limits^{n-1}_{j=0 \atop
p(j)=p(n+m)}\sum\limits_{i=j+1}^{n}{n \choose i}{j+m \choose
i-1}{i-1 \choose j}\\ =& \textstyle\sum\limits^{n-1}_{j=0 \atop
p(j)=p(n+m)}\sum\limits_{i=1}^{n-j}{n-i \choose j}{m+n-i \choose n-i}\\
=& \textstyle\sum\limits_{i=1}^{n}\sum\limits^{n-i}_{j=0 \atop
p(j)=p(n+m)}{n-i \choose j}{m+n-i \choose n-i}\\ =&
\textstyle\sum\limits_{i=1}^{n}{m+n-i \choose n-i}
\sum\limits^{n-i}_{j=0 \atop p(j)=p(n+m)}{n-i \choose j}\\
=&\textstyle\sum\limits_{i=1}^{n-1}2^{n-i-1}{m+n-i \choose n-i}+
\left\{\begin{array}{ll} 1,&\mbox{if $n+m$ is
even;}\\ 0,&\mbox{if $n+m$ is odd.}\end{array}\right.\\
=& \textstyle\sum\limits_{i=1}^{n-1}2^{i-1}{m+i \choose i}+
\left\{\begin{array}{ll} 1,&\mbox{if $n+m$ is even;}\\ 0,&\mbox{if
$n+m$ is odd.}\end{array}\right. \end{aligned}$$ where we apply
Lemma 5 and Lemma 2 in the first and the third steps respectively.
By the proof of Lemma 4, we have $\rank \, \tau^{m}+\rank \,
\tau^{m+1}= U_n^m+U_n^{m+1}+1= 2^{n-1}{n+m-1 \choose n-1}.$
\hfill$\square$

\medskip

{\bf Theorem 3.}  {\it Let} $\Lambda=kQ/I$ {\it be the exterior
algebra. Then}
$$\dim_k HH^m(\Lambda)=\left\{\begin{array}{ll}
2^n{n+m-1 \choose n-1},&\mbox{if char $k = 2$;}\\
2^{n-1}+1,&\mbox{if $m=0$, $n$ is odd and char $k \ne 2$;} \\
2^{n-1}{n+m-1 \choose n-1},& \mbox{otherwise.}
\end{array}\right.$$

\medskip

{\bf Proof.} {\it Case} char $k \ne 2$: If $m \geq 1$ then, by
Lemma 6,
$$
\begin{aligned}
\dim_k HH^m(\Lambda)=&\dim_k \, \Ker \, \tau^m-\dim_k \, \Im \, \tau^{m+1}\nonumber\\
=& \dim_k P^*_m-\dim_k \, \Im \, \tau^m-\dim_k \, \Im \, \tau^{m+1}\nonumber\\
=&\textstyle 2^n{n+m-1 \choose n-1}-\rank \, \tau^m-\rank \, \tau^{m+1}\nonumber\\
=&\textstyle 2^{n-1}{n+m-1 \choose n-1}.
\end{aligned}$$

If $m=0$ then $HH^0(\Lambda) = Z(\Lambda)$ which is the center of
$\Lambda$. Since $\Lambda$ is gradable, $\lambda \in Z(\Lambda)$
if and only if all its components belong to $Z(\Lambda)$. Thus it
is enough to consider ${\cal B} \cap Z(\Lambda)$. Note that
$x_{t_1} \cdots x_{t_i} \in Z(\Lambda)$ if and only if $x_{t_1}
\cdots x_{t_i}x_j = x_jx_{t_1} \cdots x_{t_i}$ for all $1 \leq j
\leq n$, if and only if $i=n$ or $i$ is even. Therefore
$\dim_kHH^0(\Lambda) = \dim_kZ(\Lambda)= \sum\limits^n_{i=0 \atop
{\scriptsize \mbox{$i$ even}}} {n \choose i} = 2^{n-1}$ if $n$ is
even, and $\dim_kHH^0(\Lambda) =2^{n-1}+1$ if $n$ is odd.

{\it Case} char $k=2:$ In these case all maps in the complex
$(M^{\bullet},\tau^{\bullet})$ are zero. If $m \geq 1$ then
$\dim_k HH^m(\Lambda)=\dim_k \Lambda^{{n+m-1 \choose n-1}} =
2^n{n+m-1 \choose n-1}$. If $m=0$ then
$\dim_kHH^0(\Lambda)=\dim_k\Hom_{\Lambda^e}(\Lambda^e,\Lambda)
=\dim_k \Lambda = 2^n$. \hfill{$\square$}

\medskip

{\bf Remark 1.} The case $n=2$ was obtained in \cite{BGMS}.

\medskip

{\bf Corollary 2.} {\it The Hilbert series of the exterior
algebra} $\Lambda$

$\sum\limits^{\infty}_{m=0} \dim_k HH^m(\Lambda)t^m =
\left\{\begin{array}{ll} \frac{2^{n-1}}{(1-t)^n},
&\mbox{if $n$ is even and char $k \ne 2$;} \\
\frac{2^{n-1}}{(1-t)^n}+1,&\mbox{if $n$ is odd and char $k \ne 2$;} \\
\frac{2^n}{(1-t)^n},& \mbox{if char $k = 2$.}
\end{array}\right.$

\medskip

{\bf Proof.} Note that $\sum\limits^{\infty}_{m=0} {n+m-1 \choose
n-1}t^m = \frac{1}{(1-t)^n}$. \hfill{$\square$}

\medskip

\section{Hochschild cohomology rings}

\medskip

In this section we shall determine the Hochschild cohomology rings
of the exterior algebras by generators and relations.

Now we construct another minimal projective bimodule resolution
$(P'_{\bullet},\delta'_{\bullet})$ of $\Lambda$. For each $m \geq
0$, we firstly construct elements $\{g^m_{1^{i_1} 2^{i_2}\cdots
n^{i_n}} \in \Lambda^{\otimes m} | i_1 + \cdots + i_n = m$ with
$(i_1,...,i_n) \in \mathbb{N}^n \}$: Let $g_0^0=1, g^1_1=x_1,
g^1_2=x_2, ..., g^1_n=x_n.$ Define $g^m_{1^{i_1}2^{i_2}\cdots
n^{i_n}}$ for all $m \geq 2$ inductively by
$g^m_{1^{i_1}2^{i_2}\cdots n^{i_n}} =
\sum\limits_{h=1}^ng^{m-1}_{1^{i_1}\cdots h^{i_h-1}\cdots n^{i_n}}
\otimes x_h, $ where $(i_1,...,i_n) \in \mathbb{N}^n,
i_1+i_2+\cdots +i_n=m$ and $g^{m-1}_{1^{i_1}\cdots h^{-1}\cdots
n^{i_n}}=0$ for all $1 \leq h \leq n$. It is easy to see that
$g^m_{1^{i_1}2^{i_2}\cdots n^{i_n}} = \sum\limits_{h=1}^n x_h
\otimes g^{m-1}_{1^{i_1}\cdots h^{i_h-1}\cdots n^{i_n}}$.

Let $P'_m:=\coprod\limits_{i_1 + i_2 + \cdots + i_n = m} \Lambda
\otimes g^m_{1^{i_1}2^{i_2}\cdots n^{i_n}} \otimes \Lambda
\subseteq \Lambda^{\otimes m+2}$ for $m \geq 0$, and let
$\wt{g}^m_{1^{i_1}2^{i_2}\cdots n^{i_n}}:=1 \otimes
g^m_{1^{i_1}2^{i_2}\cdots n^{i_n}} \otimes 1$ for $m \geq 1$ and
$\wt{g}^0_0=1 \otimes 1$. Note that we identify $P'_0$ with
$\Lambda \otimes \Lambda$. Define $\delta'_m: P'_m \rightarrow
P'_{m-1}$ by setting $\delta'_m(\wt{g}^m_{1^{i_1}2^{i_2}\cdots
n^{i_n}})= \sum\limits_{h=1}^n (x_h \wt{g}^{m-1}_{1^{i_1}\cdots
h^{i_h-1}\cdots n^{i_n}} +(-1)^m\wt{g}^{m-1}_{1^{i_1}\cdots
h^{i_h-1}\cdots n^{i_n}}x_h).$

\medskip

{\bf Lemma 7.} {\it The complex} $\mathbb{P}:=
(P'_{\bullet},\delta'_{\bullet}):$
$$\qquad \cdots\rightarrow P'_{m+1}
\stackrel{\delta'_{m+1}}{\longrightarrow} P'_m
\stackrel{\delta'_m}{\longrightarrow} \cdots
\stackrel{\delta'_3}{\longrightarrow} P'_2
\stackrel{\delta'_2}{\longrightarrow} P'_1
\stackrel{\delta'_1}{\longrightarrow} P'_0 \longrightarrow 0$$
{\it is a minimal projective bimodule resolution of the exterior
algebra} $\Lambda=kQ/I$.

\medskip

{\bf Proof.} Clearly the complex $\mathbb{P}=
(P'_{\bullet},\delta'_{\bullet})$ is isomorphic to the complex
$(P_{\bullet},\delta_{\bullet})$. \hfill{$\Box$}

\medskip

Applying the functor $\Hom_{\Lambda^e}(-,\Lambda)$, we have
$\mathbb{P}^*= (P'_{\bullet},\delta'_{\bullet})^* \cong
(P_{\bullet},\delta_{\bullet})^* =
(P^*_{\bullet},\delta^*_{\bullet}) \cong
(M^{\bullet},\tau^{\bullet}).$ Thus every element in $P'^*_m$ can
be represented as a linear combination of the elements in $\{\lz
e^m_{1^{i_1}2^{i_2} \cdots n^{i_n}} | \lz \in {\cal B}, i_1+i_2+
\cdots +i_n=m \} \subseteq M^m$. Throughout we do not distinguish
an element in $\Ker \, \delta'^*_{m+1} \subseteq P'^*_m$ with its
equivalent class in $HH^m(\Lambda)$.

\medskip

{\bf Lemma 8.} {\it Let} $\eta=\sum\limits_{i_1+i_2+ \cdots
+i_n=s} \lz_{i_1i_2 \cdots i_n}e^s_{1^{i_1}2^{i_2} \cdots n^{i_n}}
\in HH^s(\Lambda)$ {\it and} $\theta = \sum\limits_{j_1+j_2+
\cdots +j_n=t} \lz'_{j_1j_2 \cdots j_n}e^t_{1^{j_1}2^{j_2} \cdots
n^{j_n}} \in HH^t(\Lambda)$. {\it Then the cup product of $\eta$
and $\theta$ in $HH^{s+t}(\Lambda)$:}
$$\eta \ast \theta = \sum\limits_{i_l+j_l=h_l \atop 1 \leq l \leq
n}\lz_{i_1i_2 \cdots i_n}\lz'_{j_1j_2 \cdots
j_n}e^{s+t}_{1^{i_1+j_1}2^{i_2+j_2} \cdots n^{i_n+j_n}}.$$

\medskip

{\bf Proof.} Here we use the same strategy as that in \cite{BGMS}.
Recall that the {\it bar resolution} $\mathbb{B} =
(B_{\bullet},b_{\bullet})$ of $\Lambda$ is given by $B_m :=
\Lambda^{\otimes m+2}$ and $b_m : B_m \rightarrow B_{m-1}, \lz_0
\otimes \lz_1 \otimes \cdots \otimes \lz_{m+1} \mapsto
\sum\limits^m_{i=0} (-1)^i\lz_0 \otimes \cdots \otimes \lz_{i-1}
\otimes \lz_i\lz_{i+1} \otimes \lz_{i+2} \otimes \cdots \otimes
\lz_{m+1}$ (cf. \cite{Lod}). Define $d_i : B_m \rightarrow
B_{m-1}, \lz_0 \otimes \lz_1 \otimes \cdots \otimes \lz_{m+1}
\mapsto \lz_0 \otimes \cdots \otimes \lz_{i-1} \otimes
\lz_i\lz_{i+1} \otimes \lz_{i+2} \otimes \cdots \otimes \lz_{m+1}$
with $i=0,1,...,m$. Then $b_m = \sum\limits^m_{i=0}(-1)^id_i$.
View $P'_m$ as a submodule of $B_m$ in the natural way. The
complex $\mathbb{P}$ is a subcomplex of $\mathbb{B}$: Indeed, for
any given $\wt{g}^m_{1^{i_1}2^{i_2}\cdots n^{i_n}}$, consider the
tensors of the forms $w \otimes x_i \otimes x_i \otimes w', w
\otimes x_i \otimes x_j \otimes w'$ or $w \otimes x_j \otimes x_i
\otimes w'$ occurring in $\wt{g}^m_{1^{i_1}2^{i_2}\cdots
n^{i_n}}$, where $w$ is a $t$-fold tensor and $w'$ is an
$(m-t)$-fold tensor. Clearly, $d_t$ maps $w \otimes x_i \otimes
x_i \otimes w'$ to $0$ for all $1 \leq t \leq m-1$. If a tensor $w
\otimes x_i \otimes x_j \otimes w'$ occurs in
$\wt{g}^m_{1^{i_1}2^{i_2}\cdots n^{i_n}}$ then by definition $w
\otimes x_j \otimes x_i \otimes w'$ occurs in
$\wt{g}^m_{1^{i_1}2^{i_2}\cdots n^{i_n}}$ as well. Thus
$d_t(\wt{g}^m_{1^{i_1}2^{i_2}\cdots n^{i_n}})=0$ for all $1 \leq t
\leq m-1$. Since $g^m_{1^{i_1}2^{i_2}\cdots n^{i_n}} =
\sum\limits_{h=1}^ng^{m-1}_{1^{i_1}\cdots h^{i_h-1}\cdots n^{i_n}}
\otimes x_h = \sum\limits_{h=1}^n x_h \otimes
g^{m-1}_{1^{i_1}\cdots h^{i_h-1}\cdots n^{i_n}}$, we have
$b_m(\wt{g}^m_{1^{i_1}2^{i_2}\cdots n^{i_n}}) =
\sum\limits^m_{i=0}(-1)^id_i(\wt{g}^m_{1^{i_1}2^{i_2}\cdots
n^{i_n}}) = d_0(\wt{g}^m_{1^{i_1}2^{i_2}\cdots n^{i_n}}) + (-1)^m
d_m(\wt{g}^m_{1^{i_1}2^{i_2}\cdots n^{i_n}}) =
\delta'_m(\wt{g}^m_{1^{i_1}2^{i_2}\cdots n^{i_n}})$.

\medskip

Consider the {\it diagonal map} $\Delta : \mathbb{B} \rightarrow
\mathbb{B} \otimes _{\Lambda} \mathbb{B}$ given by $\Delta(\lz_0
\otimes \cdots \otimes \lz_{m+1}) = \sum\limits^m_{i=0}(\lz_0
\otimes \cdots \otimes \lz_i \otimes 1) \otimes_{\Lambda} (1
\otimes \lz_{i+1} \otimes \cdots \otimes \lz_{m+1}).$ The cup
product $\eta' \cup \theta'$ in the Hochschild cohomology ring
$HH^*(\Lambda)$ of two cycles $\eta'$ and $\theta'$ from
$\Hom_{\Lambda^e}(\mathbb{B},\Lambda)$ is given by the composition
$\mathbb{B} \stackrel{\Delta}{\longrightarrow} \mathbb{B}
\otimes_{\Lambda} \mathbb{B} \stackrel{\eta' \otimes
\theta'}{\longrightarrow} \Lambda \otimes _{\Lambda} \Lambda
\stackrel{\nu}{\longrightarrow} \Lambda$ where $\nu : \Lambda
\otimes _{\Lambda} \Lambda \rightarrow \Lambda$ is the
multiplication in $\Lambda$ (cf. \cite{Sa}).

\medskip

Let $\mu : \mathbb{P} \rightarrow \mathbb{B}$ be the natural
inclusion and let $\pi : \mathbb{B} \rightarrow \mathbb{P}$ be a
chain map such that $\pi \mu = 1_{\mathbb{P}}$. We have
$\Delta(\mu \mathbb{P}) \subseteq \mu \mathbb{P} \otimes
_{\Lambda} \mu \mathbb{P} \subseteq \mathbb{B} \otimes _{\Lambda}
\mathbb{B}$: Indeed, fix an $s$ with $0 \leq s \leq m$, it is easy
to see that $g^m_{1^{h_1}2^{h_2}\cdots n^{h_n}} =
\sum\limits_{i_l+j_l=h_l \atop 1 \leq l \leq n}
g^s_{1^{i_1}2^{i_2}\cdots n^{i_n}} g^{m-s}_{1^{j_1}2^{j_2} \cdots
n^{j_n}}$, here we neglect $g^0_0$ if $s=0$ or $m$. Thus we can
infer that $\Delta(\wt{g}^m_{1^{h_1}2^{h_2}\cdots n^{h_n}}) =
\sum\limits^m_{s=0}\sum\limits_{i_l+j_l=h_l \atop 1 \leq l \leq n}
\wt{g}^s_{1^{i_1}2^{i_2}\cdots n^{i_n}}
\wt{g}^{m-s}_{1^{j_1}2^{j_2}\cdots n^{j_n}}$ which implies that
$\Delta(\mu \mathbb{P}) \subseteq \mu \mathbb{P} \otimes
_{\Lambda} \mu \mathbb{P}$.

\medskip

Let $\Delta'$ be the restriction of $\Delta$ on $\mathbb{P}$. Then
$\Delta\mu = (\mu \otimes \mu)\Delta'$. Viewed as elements in
$P'^*_s$ and $P'^*_t$, $\eta$ and $\theta$ can be represented by
$\eta\pi_s$ and $\theta\pi_t$ by the bar resolution respectively.
Therefore
$$\begin{array}{lll} \eta \ast \theta &=& (\eta \ast
\theta)\pi_{s+t}\mu_{s+t}\\ &=&(\eta\pi_s \cup
\theta\pi_t)\mu_{s+t}\\ &=&\nu(\eta\pi_s \otimes
\theta\pi_t)\Delta\mu_{s+t}\\ &=&\nu(\eta\pi_s \otimes
\theta\pi_t)(\mu_s \otimes \mu_t)\Delta'\\ &=&\nu(\eta \otimes
\theta)\Delta'. \end{array}$$

\medskip

Since $\eta(\wt{g}^s_{1^{i_1}2^{i_2}\cdots n^{i_n}})= \lz_{i_1i_2
\cdots i_n}$ and $\theta(\wt{g}^t_{1^{j_1}2^{j_2}\cdots n^{j_n}})=
\lz'_{j_1j_2 \cdots j_n}$, we have $$(\eta \ast \theta)
(\wt{g}^{s+t}_{1^{h_1}2^{h_2}\cdots n^{h_n}})=
\sum\limits_{i_l+j_l=h_l \atop 1 \leq l \leq n}\lz_{i_1i_2 \cdots
i_n}\lz'_{j_1j_2 \cdots j_n}$$ for all $(h_1,...,h_n) \in
\mathbb{N}^n$ with $h_1 + \cdots + h_n=s+t$. Viewed as an element
in $M^{s+t}$, $$\eta \ast \theta = \sum\limits_{i_l+j_l=h_l \atop
1 \leq l \leq n}\lz_{i_1i_2 \cdots i_n}\lz'_{j_1j_2 \cdots
j_n}e^{s+t}_{1^{i_1+j_1}2^{i_2+j_2} \cdots n^{i_n+j_n}}.$$
\hfill{$\Box$}

\medskip

{\bf Lemma 9.} {\it Let $\Lambda=kQ/I$ be the exterior algebra and
$\char \,  k \neq 2$. Then the $k$-vector space $HH^m(\Lambda)$
has a basis $\{\lz e^m_{1^{i_1}2^{i_2} \cdots n^{i_n}} | \lz \in
{\cal B}_i, p(i)=p(m), 1 \leq i \leq n\}$.}

\medskip

{\bf Proof.} If $\lz \in {\cal B}_i$ and $p(i)=p(m)$ then
$\tau^{m+1}(\lz e^m_{1^{i_1}2^{i_2}\cdots n^{i_n}})=0$. Thus we
have $\{\lz e^m_{1^{i_1}2^{i_2} \cdots n^{i_n}} | \lz \in {\cal
B}_i, p(i)=p(m), 1 \leq i \leq n\} \subseteq \Ker \, \tau^{m+1}$.
Clearly $\Im \, \tau^m$ is contained in the subspace of $\Ker \,
\tau^{m+1}$ generated by $\{\lz e^m_{1^{i_1}2^{i_2} \cdots
n^{i_n}} | \lz \in {\cal B}_i, p(i) \neq p(m), 1 \leq i \leq n\}$.
By Theorem 3, we have $\dim_k HH^m(\Lambda) = 2^{n-1} {n+m-1
\choose n-1}$. Thus the Lemma holds. \hfill{$\Box$}

\medskip

It follows from Lemma 9 that, as a $k$-algebra, $HH^*(\Lambda)$ is
generated by $\{x_ix_j|1 \leq i < j \leq n\} \cup \{x_pe^1_q|1
\leq p,q \leq n\} \cup \{e^2_{st}|1 \leq s \leq t \leq n\}$ which
satisfy all relations in the following Table H.

\newpage

\bcen{\bf Table H}\ecen

\medskip

$\begin{array}{lll}
(H1.1)&(x_ix_j)(x_sx_t)=(x_sx_t)(x_ix_j)&\mbox{if $i<j$ and $s<t$}\\
(H1.2)&(x_ix_j)(x_sx_t)=0&\mbox{if $\{i,j\}\cap\{s,t\} \neq \emptyset$}\\
(H1.3)&(x_ix_j)(x_sx_t)=-(x_ix_s)(x_jx_t)&\mbox{if $i<s<j<t$}\\
(H1.4)&(x_ix_j)(x_sx_t)=(x_ix_s)(x_tx_j)&\mbox{if $i<s<t<j$}\\
(H1.5)&(x_ix_j)(x_sx_t)=-(x_sx_i)(x_tx_j)&\mbox{if $s<i<t<j$}\\
(H1.6)&(x_ix_j)(x_sx_t)=(x_sx_i)(x_jx_t)&\mbox{if $s<i<j<t$}\\
&&\\
(H2.1)&(x_ix_j)(x_se^1_t)=(x_se^1_t)(x_ix_j)&\mbox{if $i<j$}\\
(H2.2)&(x_ix_j)(x_se^1_t)=0&\mbox{if $s \in \{i,j\}$}\\
(H2.3)&(x_ix_j)(x_se^1_t)=(x_sx_i)(x_je^1_t)&\mbox{if $s<i<j$}\\
(H2.4)&(x_ix_j)(x_se^1_t)=-(x_ix_s)(x_je^1_t)&\mbox{if $i<s<j$}\\
&&\\
(H3.1)&(x_ix_j)e^2_{st}=e^2_{st}(x_ix_j)&\mbox{if $i<j$ and $s \leq t$}\\
&&\\
(H4.1)&(x_ie^1_j)(x_se^1_t)=0&\mbox{if $i=s$}\\
(H4.2)&(x_ie^1_j)(x_se^1_t)=(x_ix_s)e^2_{jt}&\mbox{if $i<s$ and $j \leq t$}\\
(H4.3)&(x_ie^1_j)(x_se^1_t)=(x_ix_s)e^2_{tj}&\mbox{if $i<s$ and $t \leq j$}\\
(H4.4)&(x_ie^1_j)(x_se^1_t)=-(x_sx_i)e^2_{jt}&\mbox{if $s<i$ and $j \leq t$}\\
(H4.5)&(x_ie^1_j)(x_se^1_t)=-(x_sx_i)e^2_{tj}&\mbox{if $s<i$ and $t \leq j$}\\
&&\\
(H5.1)&(x_ie^1_j)e^2_{st}=e^2_{st}(x_ie^1_j)&\mbox{if $s \leq t$}\\
(H5.2)&(x_ie^1_j)e^2_{st}=(x_ie^1_s)e^2_{jt}&\mbox{if $s<j \leq t$}\\
(H5.3)&(x_ie^1_j)e^2_{st}=(x_ie^1_s)e^2_{tj}&\mbox{if $s<t \leq j$}\\
&&\\
(H6.1)&e^2_{ij}e^2_{st}=e^2_{st}e^2_{ij}&\mbox{if $i \leq j$ and $s \leq t$}\\
(H6.2)&e^2_{ij}e^2_{st}=e^2_{is}e^2_{jt}&\mbox{if $i \leq s \leq j \leq t$}\\
(H6.3)&e^2_{ij}e^2_{st}=e^2_{is}e^2_{tj}&\mbox{if $i \leq s \leq t \leq j$}\\
(H6.4)&e^2_{ij}e^2_{st}=e^2_{si}e^2_{tj}&\mbox{if $s \leq i \leq t \leq j$}\\
(H6.5)&e^2_{ij}e^2_{st}=e^2_{si}e^2_{jt}&\mbox{if $s \leq i \leq j
\leq t$}
\end{array}$

\medskip

Replace $x_ix_j, x_pe^1_q,e^2_{st}$ with $u_{ij}, v_{pq},w_{st}$
respectively, we have the following Table F.

\newpage

\bcen{\bf Table F}\ecen

\medskip

$\begin{array}{lll}
(F1.1)&u_{ij}u_{st}=u_{st}u_{ij}&\mbox{if $i<j$ and $s<t$}\\
(F1.2)&u_{ij}u_{st}=0&\mbox{if $\{i,j\}\cap\{s,t\} \neq \emptyset$}\\
(F1.3)&u_{ij}u_{st}=-u_{is}u_{jt}&\mbox{if $i<s<j<t$}\\
(F1.4)&u_{ij}u_{st}=u_{is}u_{tj}&\mbox{if $i<s<t<j$}\\
(F1.5)&u_{ij}u_{st}=-u_{si}u_{tj}&\mbox{if $s<i<t<j$}\\
(F1.6)&u_{ij}u_{st}=u_{si}u_{jt}&\mbox{if $s<i<j<t$}\\
&&\\
(F2.1)&u_{ij}v_{st}=v_{st}u_{ij}&\mbox{if $i<j$}\\
(F2.2)&u_{ij}v_{st}=0&\mbox{if $s \in \{i,j\}$}\\
(F2.3)&u_{ij}v_{st}=u_{si}v_{jt}&\mbox{if $s<i<j$}\\
(F2.4)&u_{ij}v_{st}=-u_{is}v_{jt}&\mbox{if $i<s<j$}\\
&&\\
(F3.1)&u_{ij}w_{st}=w_{st}u_{ij}&\mbox{if $i<j$ and $s \leq t$}\\
&&\\
(F4.1)&v_{ij}v_{st}=0&\mbox{if $i=s$}\\
(F4.2)&v_{ij}v_{st}=u_{is}w_{jt}&\mbox{if $i<s$ and $j \leq t$}\\
(F4.3)&v_{ij}v_{st}=u_{is}w_{tj}&\mbox{if $i<s$ and $t \leq j$}\\
(F4.4)&v_{ij}v_{st}=-u_{si}w_{jt}&\mbox{if $s<i$ and $j \leq t$}\\
(F4.5)&v_{ij}v_{st}=-u_{si}w_{tj}&\mbox{if $s<i$ and $t \leq j$}\\
&&\\
(F5.1)&v_{ij}w_{st}=w_{st}v_{ij}&\mbox{if $s \leq t$}\\
(F5.2)&v_{ij}w_{st}=v_{is}w_{jt}&\mbox{if $s<j \leq t$}\\
(F5.3)&v_{ij}w_{st}=v_{is}w_{tj}&\mbox{if $s<t \leq j$}\\
&&\\
(F6.1)&w_{ij}w_{st}=w_{st}w_{ij}&\mbox{if $i \leq j$ and $s \leq t$}\\
(F6.2)&w_{ij}w_{st}=w_{is}w_{jt}&\mbox{if $i \leq s \leq j \leq t$}\\
(F6.3)&w_{ij}w_{st}=w_{is}w_{tj}&\mbox{if $i \leq s \leq t \leq j$}\\
(F6.4)&w_{ij}w_{st}=w_{si}w_{tj}&\mbox{if $s \leq i \leq t \leq j$}\\
(F6.5)&w_{ij}w_{st}=w_{si}w_{jt}&\mbox{if $s \leq i \leq j \leq
t$}
\end{array}$

\medskip

{\bf Theorem 4.} {\it Let $Q'$ be the quiver with one vertex $1$
and $2n^2$ loops $\{u_{ij}|1 \leq i < j \leq n\} \cup \{v_{pq}|1
\leq p,q \leq n\} \cup \{w_{st}|1 \leq s \leq t \leq n\}$. Let
$I'$ be the ideal of $kQ'$ generated by all relations in Table F.
Then the Hochschild cohomology ring of the exterior algebra}
$\Lambda=kQ/I$, $$HH^*(\Lambda) =
\left\{\begin{array}{ll} kQ'/I',&\mbox{if $\char \, k \neq 2$;}\\
\Lambda[z_1,....,z_n],&\mbox{if $\char \, k = 2$.}
\end{array}\right.$$

\medskip

{\bf Proof.} {\it Case} $\char \, k \neq 2$: Firstly, as a
$k$-algebra $HH^*(\Lambda)$ is generated by $\{x_ix_j|1 \leq i < j
\leq n\} \cup \{x_pe^1_q|1 \leq p,q \leq n\} \cup \{e^2_{st}|1
\leq s \leq t \leq n\}$. Thus we have an epimorphism of
$k$-algebras $\psi: kQ' \rightarrow HH^*(\Lambda)$ which maps
$u_{ij},v_{pq},w_{st}$ to $x_ix_j,x_pe^1_q,e^2_{st}$ respectively.
Comparing the relations in Table H with those in Table F, we have
$I' \subseteq \Ker \psi$. Hence $\psi$ induces an epimorphism
$\varphi : kQ'/I' \rightarrow HH^*(\Lambda)$.

Apply the relations in Table H one by one, it is not difficult to
see that the set $\{(x_{i_1}x_{i_2}) \cdots
(x_{i_{2l-1}}x_{i_{2l}})e^2_{j_1j_2} \cdots
e^2_{j_{2r-1}j_{2r}}|1<i_1< i_2<\cdots<i_{2l-1}<i_{2l}, j_1 \leq
j_2 \leq \cdots \leq j_{2r}\} \cup \{(x_{i_1}x_{i_2}) \cdots
(x_{i_{2l-1}}x_{i_{2l}})(x_{2l+1}e^1_{j_1})e^2_{j_2j_3} \cdots
e^2_{j_{2r}j_{2r+1}} | 1<i_1<i_2<\cdots<i_{2l-1}<i_{2l}<i_{2l+1},
j_1 \leq j_2 \leq \cdots \leq j_{2r} \leq j_{2r+1}\}$ is a
$k$-basis of $HH^*(\Lambda)$.

Apply the relations in Table F one by one, it is not difficult to
see that the set $\{u_{i_1i_2} \cdots u_{i_{2l-1}i_{2l}}w_{j_1j_2}
\cdots w_{j_{2r-1}j_{2r}}| 1<i_1< i_2<\cdots<i_{2l-1}<i_{2l}, j_1
\leq j_2 \leq \cdots \leq j_{2r}\} \cup \{u_{i_1i_2} \cdots
u_{i_{2l-1}i_{2l}}v_{(2l+1)j_1}w_{j_2j_3} \cdots
w_{j_{2r}j_{2r+1}} | 1<i_1< i_2<\cdots<i_{2l-1}<i_{2l}<i_{2l+1},
j_1 \leq j_2 \leq \cdots \leq j_{2r} \leq j_{2r+1}\}$ is a
$k$-basis of $kQ'/I'$.

Since $\varphi$ maps a basis element to a basis element,
$\varphi$ is also a monomorphism. Hence it is an isomorphism.

\medskip

{\it Case} $\char \, k =2$: The differential in the complex
$\Hom_{\Lambda^e}(\mathbb{P},\Lambda)$ is $0$. All maps in
$\Hom_{\Lambda^e}(\mathbb{P},\Lambda)$ represent nonzero elements
in $HH^*(\Lambda)$. Let $\phi : \Lambda[z_1,....,z_n] \rightarrow
HH^*(\Lambda), z_i \mapsto e^1_i$. Then $\phi$ is surjective and
thus it is injective restricted to each degree. Hence it is an
isomorphism. \hfill{$\Box$}

\medskip

{\bf Remark 2.} The case $n=2$ was obtained in \cite{BGMS}.

\end{document}